\documentclass[hidelinks]{article}

\usepackage{graphicx, bm, amsmath, siunitx, longtable, tabularx, float, mathtools, booktabs, eqnarray, color, xcolor, soul, amsfonts, amssymb, appendix, tensor, subcaption, enumitem, makecell, listings, dsfont, lscape, diagbox, attachfile}
\usepackage[square,numbers]{natbib}
\usepackage[justification=centering]{caption}

\newcommand{\tfc}{Theory of Functional Connections}
\newcommand{\ces}{constrained expressions}
\newcommand{\ce}{constrained expression}
\newcommand{\p}[2]{\prescript{(#1)}{}{#2}}

\newcommand{\C}[1]{\tensor*[^{}_{}]{\mathfrak{C}}{^{}_{#1}}}
\newcommand{\pC}[2]{\tensor*[^{(#1)}_{}]{\mathfrak{C}}{^{}_{#2}}}

\newcommand{\andd}{\quad \text{and} \quad}
\newcommand{\R}{\mathbb{R}}
\newcommand{\F}{\mathbb{F}}

\title{The \tfc\ on Vector Spaces}
\author{Carl Leake}
\date{}

\begin{document}

\maketitle
\begin{abstract}
\noindent The \tfc\ (TFC) is most often used for constraints over the field of real numbers. However, previous works have shown that it actually extends to arbitrary fields. The evidence for these claims is restricting oneself to the field of real numbers is unnecessary because all of the theorems, proofs, etc. for TFC apply as written to fields in general. Here, that notion is taken a step further, as fields themselves are unnecessarily restrictive. The theorems, proofs, etc. of TFC apply to vector spaces, not just fields. Moreover, the inputs of the functions/constraints do not even need to be restricted to vector spaces; they can quite literally be anything. In this note, the author shows the more general inputs/outputs of the various symbols that comprise TFC. Examples with accompanying, open-access Python code are included to aid the reader's understanding.
\end{abstract}

\section{Introduction}
The main idea and mathematical justification given here is simple. Hence, this is written as a note rather than an article, as the author does not feel a full article is warranted.

Proofs, theorems, definitions, examples, properties, etc. of the \tfc\ (TFC) (hereafter referred to as ``the results of TFC'') are typically done using the field of real numbers since they are used the most in its applications. However, as conjectured in Appendix I of Ref. \cite{LeakeDissertation} and Appendix B of Ref. \cite{TFC_Book}, the results of TFC apply to any field\footnote{These references also include examples over the field of complex numbers and finite fields. In addition, the reference documentation of the TFC Python module provides an example of a complex ODE solved via TFC \cite{TfcGithub}.}. The evidence for this conjecture is TFC only uses properties and axioms defined for all fields. In other words, restricting oneself to the field of real numbers is unnecessary, as none of the results of TFC require properties or axioms specific to the field of real numbers.

This note shows that fields are still unnecessarily restrictive. The only time that a multiplicative inverse is required is when inverting the support matrix. All other operations only require the properties of vector spaces. This means that the results of TFC as already defined in previous works actually apply to functions whose outputs are vector spaces defined over a field. Notice that only the outputs of these functions need to be vector spaces. Since the inputs of the functions are not used in the constrained expression, they have no restriction, at least not in the univariate case. 

For example, consider the algorithm used to generate a univariate constrained expression.
Given a set of constraints
\begin{equation*}
    \C{i}[u(x)] = \kappa_i, \quad i\in 1,\dots,k
\end{equation*}
the univariate constrained expression that satisfies these constraints can be derived using:
\begin{enumerate}
    \item Derive the projection functions
    \begin{equation*}
        \rho_i(x,g(x)) = \kappa_i - \C{i}[g(x)].
    \end{equation*}
    \item Choose $k$ support functions, $s_i(x)$.
    \item Derive the switching functions
    \begin{equation*}
        \phi_i(x) = s_j(x)\alpha_{ji} \quad \text{where} \quad \alpha_{ji} = \Big(\C{j}[s_i(x)]\Big)^{-1}
    \end{equation*}
    \item Create the constrained expression
    \begin{equation*}
        u(x,g(x)) = g(x) + \phi_i(x) \rho_i(x,g(x)).
    \end{equation*}
\end{enumerate}

In this algorithm, a multiplicative inverse is only required when solving for $\alpha_{ij}$. All other operations only require the axioms and properties of a vector space, e.g., addition, an additive inverse, distribution of a scalar over addition, etc.

\section{Univariate TFC}
In this section, the input and output types for the symbols associated with univariate constrained expressions are shown. 

Let $V$ be a vector space defined over a field $\F$. Moreover, let $A$ denote any type. It may be easiest to think of ``any'' here in the same way as the ``Any'' type is defined in Python. This can truly be anything: a number, a string, a picture, etc. Using this definition, one can define the input and output types of all the symbols used in the above algorithm. 
\begin{align*}
\C{i}\colon& (A\mapsto V) \mapsto V\\
\kappa\colon& V\\
g(x)\colon& A\mapsto V\\
s_i(x)\colon& A\mapsto \F\\
\phi_i(x)\colon& A\mapsto \F\\
\rho_i(x,g(x))\colon& (A, A\mapsto V) \mapsto V
\end{align*}
In the above, the notation $(A, A\mapsto V) \mapsto V$ indicates that the inputs are $A$ and a function of type $A\mapsto V$ and the output is $V$.

\section{Multivariate TFC}
The same idea applies to multivariate TFC. The only added restriction is the input types need to be orthogonal in the sense that the univariate \ce\ from one variable should not interfere with the univariate \ce\ from another.

Let $V$ be a vector space defined over the field $\F$. Moreover, let $A_j$ denote types that are orthogonal where $j\in 1\dots n$ and $n$ is the number of dimensions. Then, for  the $k$-th dimension:
\begin{align*}
\pC{k}{i}\colon& \Big((A_1, \dots, A_n) \mapsto V\Big) \mapsto \Big((A_1, \dots, A_{k-1}, A_{k+1}, \dots, A_n) \mapsto V \Big)\\
\kappa\colon& (A_1, \dots, A_{k-1}, A_{k+1}, \dots, A_n) \mapsto V\\
g(x)\colon& (A_1, \dots, A_n) \mapsto V\\
s_i(x)\colon& (A_1, \dots, A_n) \mapsto \F\\
\phi_i(x)\colon& (A_1, \dots, A_n) \mapsto \F\\
\rho_i(x,g(x))\colon& \Big((A_1, \dots, A_n), (A_1, \dots, A_n)\mapsto V\Big) \mapsto V
\end{align*}

\section{Examples}
This section includes examples that TFC applied on various vector spaces. Most of the examples here are done using univariate TFC, since that makes the examples simpler, shorter, and easier to follow. However, one example is included using multivariate TFC for completeness. All examples have accompanying, open-access Python code, which can be found in Ref. \cite{tfc_vector_spaces_github}.

\subsection{Univariate \texorpdfstring{$\R^{2\times 3} \mapsto \R^2$}{R2x3 to R2} }
Let $x \in \mathbb{R}^{2\times 3}$ and the output of $u(x) \in R^2$. Consider the constraints
\begin{align*}
    u\Bigg(\begin{bmatrix} 1 & 3 & 5 \\ 2 & 4 & 7 \end{bmatrix}\Bigg) = \begin{Bmatrix} 3 \\ 4 \end{Bmatrix} \andd u\Bigg(\begin{bmatrix} 1 & 0 & 1 \\ 0 & 1 & 0 \end{bmatrix}\Bigg) = u\Bigg(\begin{bmatrix} 0 & 0 & 1 \\ 1 & 0 & 0 \end{bmatrix}\Bigg)
\end{align*}

Let the support functions be
\begin{align*}
    s_1(x) = \begin{Bmatrix} 1 & 0 \end{Bmatrix} x \begin{Bmatrix} 1 \\ 0 \\ 0 \end{Bmatrix} \andd s_2(x) = \begin{Bmatrix} 1 & 0 \end{Bmatrix} x \begin{Bmatrix} 0 \\ 1 \\ 0 \end{Bmatrix}.
\end{align*}
Then, following the standard TFC \ce\ algorithm yields:
\begin{align*}
    \alpha_{ij} &= \Big(\C{i}[s_j(x)]\Big)^{-1} = \begin{bmatrix} 1 & 3 \\ 1 & 0 \end{bmatrix}^{-1}\\
    \alpha_{ij} &= \begin{bmatrix} 0 & 1 \\ \frac{1}{3} & -\frac{1}{3} \end{bmatrix}
\end{align*}
\begin{align*}
   u(x, g(x)) =&\ g(x) + \frac{1}{3} \Bigg( \begin{Bmatrix} 1 & 0 \end{Bmatrix} x \begin{Bmatrix} 0 \\ 1 \\ 0 \end{Bmatrix} \Bigg) \Bigg( \begin{Bmatrix} 3 \\ 4 \end{Bmatrix} - g\Bigg(\begin{bmatrix} 1 & 3 & 5 \\ 2 & 4 & 7 \end{bmatrix}\Bigg)\Bigg) \\ 
   &+ \Bigg(\begin{Bmatrix} 1 & 0 \end{Bmatrix} x \begin{Bmatrix} 1 \\ 0 \\ 0 \end{Bmatrix} - \frac{1}{3} \begin{Bmatrix} 1 & 0 \end{Bmatrix} x \begin{Bmatrix} 0 \\ 1 \\ 0 \end{Bmatrix} \Bigg) \Bigg(g\Bigg(\begin{bmatrix} 0 & 0 & 1 \\ 1 & 0 & 0 \end{bmatrix}\Bigg) - g\Bigg(\begin{bmatrix} 1 & 0 & 1 \\ 0 & 1 & 0 \end{bmatrix}\Bigg)\Bigg)
\end{align*}

This \ce\ satisfies the constraints. This can easily be verified by evaluating the \ce\ at the constraint points. Doing so is simple enough in this example, so it is left as an exercise for the reader.

\subsection{Univariate \texorpdfstring{string $\mapsto \mathbb{C}^2$}{string to C2} }
For this example, let the input $x$ be a member of the set of all possible strings whose characters are the letters a through z (lowercase) and the output of $u(x)\in \mathbb{C}^2$. An example of such a function would be a cipher that maps each character in the string to a complex number and takes the sum; two different ciphers could be used to map the input string to an element of $\mathbb{C}^2$. 

Consider the constraints
\begin{align*}
    u(``tfc") = \begin{Bmatrix} 42i \\ 19 + 79i \end{Bmatrix} \andd u(``dont") = u(``panic")
\end{align*}

Let the first support function be the constant function
\begin{align*}
    s_1(x) = 4 + 2i.
\end{align*}
Let the second support function use two ciphers to map the characters of the input string to numbers; the outputs of each cipher make up the real and imaginary portion of the output. The first cipher maps a to 1, b to 2, and so on, and the second maps them to $i^k$ where $k$ is the number from the first cipher. For example, the first cipher would map ``tfc'' to $(t = 20) + (f = 6) + (c = 3) = 29$ and the second cipher would map ``tfc'' to $i^{20} + i^6 + i^3 = -i$. Hence, 
\begin{align*}
    s_2(``tfc") = 29 -i.
\end{align*}

Following the standard TFC \ce\ algorithm yields
\begin{align*}
    \alpha_{ij} &= \Big(\C{i}[s_j(x)]\Big)^{-1} = \begin{bmatrix} 4+2i & 29-i \\ 0 & 11-2i \end{bmatrix}^{-1}\\
    \alpha_{ij} &= \begin{bmatrix} 0.2-0.1i & -0.5512 + 0.1816i \\ 0 & 0.088 + 0.016i \end{bmatrix}
\end{align*}
\begin{align*}
   u(x, g(x)) =&\ g(x) + (0.2 - 0.1i) (4 + 2i) \Bigg(\begin{Bmatrix} 42i \\ 19 + 79i \end{Bmatrix} - g(``tfc")\Bigg)\\
   &+\Big((-0.5512 + 0.1816i) (4 + 2i) + (0.088 + 0.016i)s_2(x)\Big) \Big(g(``panic") - g(``dont")\Big)
\end{align*}

Simplifying a bit yields,
\begin{align*}
   u(x, g(x)) =&\ g(x) + \begin{Bmatrix} 42i \\ 19 + 79i \end{Bmatrix} - g(``tfc")\\
   &+\Big(-2.568 - 0.376i + (0.088 + 0.016i)s_2(x)\Big) \Big(g(``panic") - g(``dont")\Big)
\end{align*}

This constrained expression satisfies the constraints as expected.
\begin{align*}
    u(``tfc", g(x)) =&\ g(``tfc'') + \begin{Bmatrix} 42i \\ 19 + 79i \end{Bmatrix} - g(``tfc")\\
    &\quad +\Big(-2.568 - 0.376i + (0.088 + 0.016i)(29-i)\Big) \Big(g(``panic") - g(``dont")\Big)\\
    =&\ \begin{Bmatrix} 42i \\ 19 + 79i \end{Bmatrix}+\Big(-2.568 - 0.376i + 2.568 + 0.376i\Big)\Big(g(``panic") - g(``dont")\Big)\\
    =&\ \begin{Bmatrix} 42i \\ 19 + 79i \end{Bmatrix}
\end{align*}
\begin{align*}
    u(``dont", g(x)) &- u(``panic", g(x)) =\\  
    &\ g(``dont'') - g(``panic'')\\
    &+\Big(-2.568 - 0.376i + (0.088 + 0.016i)(54-i)\Big) \Big(g(``panic") - g(``dont")\Big)\\
    &-\Big(-2.568 - 0.376i + (0.088 + 0.016i)(43+i)\Big) \Big(g(``panic") - g(``dont")\Big)\\
    =&\ g(``dont'') - g(``panic'')\\
    &+\Big(2.2 + 0.4i )\Big) \Big(g(``panic") - g(``dont")\Big)-\Big(1.2 + 0.4i)\Big) \Big(g(``panic") - g(``dont")\Big)\\
    =&\ g(``dont'') - g(``panic'') + g(``panic'') - g(``dont'')\\
    =&\ 0 \\
    \Rightarrow u(``dont", g(x)) &= u(``panic", g(x))
\end{align*}

\subsection{Univariate \texorpdfstring{$\F \mapsto (\F \mapsto \F$)}{F to (F to F)}}
Of course, functions that go from a fixed set to a field form a vector space. Hence, TFC applies here as well. In this example, let $x\in \F$ where $\F$ is a finite field containing the four elements $\{0, 1, A, B\}$. The addition and multiplication tables for this field are shown in Tables \ref{tab:FiniteFieldPlus} and \ref{tab:FiniteFieldMult}. Furthermore, let $u(x) \in (\F \mapsto \F)$. In other words, the inputs are from the given finite field and the outputs are functions that map values on that finite field to other values on that finite field. To easily distinguish the outputs of $u(x)$, let these functions have input $t\in\F$. For example, a constraint might look like $u(1) = t + A$.

\begin{table}[H]
    \centering
    \caption{Addition table.}
    \label{tab:FiniteFieldPlus}
    \begin{tabular}{|c|c|c|c|c|}
        \hline
        $+$ & $0$ & $1$ & $A$ & $B$ \\\hline
        $0$ & $0$ & $1$ & $A$ & $B$ \\\hline
        $1$ & $1$ & $0$ & $B$ & $A$ \\\hline
        $A$ & $A$ & $B$ & $0$ & $1$ \\\hline
        $B$ & $B$ & $A$ & $1$ & $0$ \\\hline
    \end{tabular}
\end{table}
\begin{table}[H]
    \centering
    \caption{Multiplication table.}
    \label{tab:FiniteFieldMult}
    \begin{tabular}{|c|c|c|c|c|}
        \hline
         $*$ & $0$ & $1$ & $A$ & $B$ \\\hline
         $0$ & $0$ & $0$ & $0$ & $0$ \\\hline
         $1$ & $0$ & $1$ & $A$ & $B$ \\\hline
         $A$ & $0$ & $A$ & $B$ & $1$ \\\hline
         $B$ & $0$ & $B$ & $1$ & $A$ \\\hline
    \end{tabular}
\end{table}

Consider the constraints
\begin{equation*}
    u(1) = t + A \andd u(A) - u(B) = t
\end{equation*}

Let the support functions be
\begin{equation*}
    s_1(x) = 1 \andd s_2(x) = x.
\end{equation*}
Then, the constrained expression is derived as
\begin{align*}
    \alpha_{ij} &= \Big(\C{i}[s_j(x)]\Big)^{-1} = \begin{bmatrix} 1 & 1 \\ 0 & 1 \end{bmatrix}^{-1}\\
    \alpha_{ij} &= \begin{bmatrix} 1 & 1 \\ 0 & 1 \end{bmatrix}
\end{align*}
\begin{align*}
    u(x,g(x)) = g(x) + 1\Big(t + A - g(1)\Big) + (1 + x)\Big(t + g(B) - g(A)\Big)
\end{align*}
As expected, the \ce\ satisfies the constraints regardless of $g(x)$. To verify, let's substitute the constraints and simplify
\begin{align*}
    u(1,g(x)) &= g(1) + 1\Big(t+A-g(1)\Big) + (0)\Big(t + g(B) - g(A)\Big)\\
    &= g(1) + t + A - g(1)\\
    &= t + A
\end{align*}
\begin{align*}
    u(A, g(x)) - u(B, g(x)) &= g(A) - g(B) + (1 - 1)\Big(t+A-g(1)\Big) + (1 + A + 1 - B)\Big(t + g(B) - g(A)\Big)\\
    &= g(A) - g(B) + (1)\Big(t - g(B) - g(A)\Big)\\
    &= g(A) - g(B) + t -g(B) - g(A)\\
    &= t
\end{align*}

\subsection{Multivariate \texorpdfstring{$(\R, \R) \mapsto (\R \mapsto \R)$}{(R, R) to (R to R)}}

Consider the set of constraints
\begin{equation*}
    u(0,y) = \sin(t), \quad u(1,y) + u(2,y) = u(5,y) + u(4,y), \andd u(x,1) = \sin(t)
\end{equation*}

For the \ce\ on $x$:
\begin{equation*}
    s_1(x) = 1 \andd s_2(x) = x
\end{equation*}
\begin{align*}
    \alpha_{ij} &= \Big(\pC{1}{i}[s_j(x)]\Big)^{-1} = \begin{bmatrix} 1 & 0 \\ 0 & -6 \end{bmatrix}^{-1}\\
    \alpha_{ij} &= \begin{bmatrix}1 & 0 \\ 0 & - \frac{1}{6}\end{bmatrix}
\end{align*}
\begin{align*}
    \p{1}{u}(x,y,g(x,y)) = g(x,y) + \Big(\sin(t) - g(0,y)\Big) - \frac{x}{6} \Big(g(4,y) + g(5,y) - g(1,y) - g(2,y)\Big).
\end{align*}
For the \ce\ on $y$:
\begin{equation*}
    s_1(y) = 1
\end{equation*}
\begin{align*}
    \alpha_{ij} &= \Big(\pC{2}{i}[s_j(y)]\Big)^{-1} = \begin{bmatrix} 1 \end{bmatrix}^{-1}\\
    \alpha_{ij} &= \begin{bmatrix}1\end{bmatrix}
\end{align*}
\begin{align*}
    \p{2}{u}(x,y,g(x,y)) = g(x,y) + \sin(t) - g(x,1) 
\end{align*}

The full constrained expression is given by
\begin{equation*}
    u(x,y,g(x,y)) = \p{2}{u}(x,y,\p{1}{u}(x,y,g(x,y)))
\end{equation*}
Although multivariate \ces\ are not typically written out completely (it is unnecessary since the recursive form can be used in code), it is done below for completeness. 
\begin{align*}
    u(x,y,g(x,y)) =&\ g(x,y) + \sin(t) - g(x,1) - g(0,y) + g(0,1) \\
    & + \frac{x}{6} \Bigg(\Big(g(4,1) + g(5,1) - g(1,1) - g(2,1)\Big) - \Big(g(4,y) + g(5,y) - g(1,y) - g(2,y)\Big)\Bigg) 
\end{align*}

This \ce\ satisfies the constraints regardless of how $g(x)$ is chosen. One can verify this by evaluating it at the constraint points. Doing so is simple enough in this example, so it is left as an exercise for the reader.

\bibliographystyle{unsrtnat}
\bibliography{refs}

\end{document}